\documentclass{article}

\usepackage{amsfonts,amsmath, amsthm, amscd, amssymb, graphicx, color, mathrsfs, stmaryrd}
\usepackage{url}

\begin{document}

\newtheorem{thm}{Theorem}[section]
\newtheorem{theo}[thm]{Theorem}
\newtheorem{prop}[thm]{Proposition}
\newtheorem{prty}[thm]{Property}
\newtheorem{coro}[thm]{Corollary}
\newtheorem{lema}[thm]{Lemma}
\newtheorem{defi}[thm]{Definition}
\newtheorem{ejem}[thm]{Example}
\newtheorem{rema}[thm]{Remark}
\newtheorem{fact}[thm]{Fact}
\newtheorem{open}[thm]{Problem}

\title{A covering index for Banach spaces}
\author{M. Raja\thanks{This research has been supported by:  Fundaci\'on S\'eneca -- ACyT Regi\'on de Murcia, project 21955/PI/22; and grant PID2021-122126NB-C32  funded by MCIN/AEI/ 10.13039/501100011033 and by “ERDF A way of making Europe”.}}
\date{November, 2025}

\maketitle

\begin{abstract}
We introduce a new isomorphic quantity for Banach spaces, the index $\Theta_X$, based on finite convex coverings of the unit ball. 
This index is closely related to the asymptotic moduli of uniform convexity and uniform smoothness, so that it can be calculated for several classical Banach spaces.
\end{abstract}

\section{Introduction}

Since the foundation of Banach space theory, it had been thought that every infinite-dimensional space contained some copy of the classical $c_0$ or $\ell_p$ spaces, however the explicit question can only be traced to the 70's, see \cite{LT} for instance. Tsirelson example proved that it is not the case, however the search for classical spaces inside general infinite dimensional spaces still went on, somehow, in weakened fashions such as type and cotype, or spreading models (Krivine's theorem). These notions involve sequences, finite or infinite, of vectors.
There exists sequence-free quantitative methods to explore the geometry of a Banach space and compare it with that of $c_0$ or $\ell_p$, for instance, the Szlenk-type indices, but the definitions are quite technical, see \cite{Lancien} for instance.\\

In this note we will propose a simpler notion, maybe even more natural, that will allow to compare an arbitrary real Banach space $X$ to the classical spaces $c_0$ or $\ell_p$.  Let $B_X$ denote the unit ball of $X$. 
The idea is to measure how may shrink, with respect to a suitable measure, the pieces of a finite cover of $B_X$ by convex sets as its number increases. That will lead to our covering index $\Theta_{X}$.
We need to introduce a few definitions.
The 
{\it essential inradius} $\varrho(A)$ of a set $A \subset X$ is the number
$$  \varrho(A) = \sup \{r>0: \exists x \in A, \, \exists Y \subset X, \, {\mathcal D}(X/Y) < \infty, ~x+rB_Y \subset A \} ,$$
where $B_Y = B_X \cap Y$ and ${\mathcal D}$ stands for dimension, so we are measuring radii of finitely codimensional balls.
The essential inradius should not be confused with the {\it asymptotic inradius} that is computed similarly but running on finite-dimensional subspaces of arbitrarily high dimension, see \cite{raja2}. Let ${\mathcal C}(X)$ denote the family of bounded closed convex subsets of $X$. Consider the following covering index 
$$ \Theta_X( n) = \inf \Big\{ \max_{1 \leq k \leq n}  \varrho(A_k) : (A_k)_{k=1}^n \subset {\mathcal C}(X), 
\bigcup_{k=1}^n A_k= B_X \Big\}. $$
That means, for every convex decomposition of $B_X$ into $n$ sets, there is one of them that contains infinite codimensional balls of radius arbitrarily close to $ \Theta(n)$, if not bigger.\\

Evidently, for $X$ an infinite-dimensional Banach space 
$(\Theta_X( n))_{n=1}^\infty \subset (0,1]$ is a non-increasing  sequence. Moreover, the information provided by 
$\Theta_X$ is of isomorphic nature. Indeed, it is not difficult to prove that if $\tilde{X}$ denotes a renorming of $X$ with a $\lambda$-equivalent norm, then
$$ \lambda^{-2} \Theta_X( n)  \leq \Theta_{\tilde{X}}(n) \leq  \lambda^{2} \Theta_X( n) . $$
Another useful observation is that in the definition of $ \Theta_X( n) $ we can only consider coverings where all the sets $A_k$ have nonempty interior, \cite[Lemma 2.4]{raja2}. That would allow the use of perturbation arguments to impose some restrictions to the finite codimensional subspaces $Y$ appearing in the definition of $\varrho(A_k)$.\\

Two examples given in \cite{raja2} can be reformulated as follows:
\begin{itemize}
\item[(a)] $\Theta_{c_0}(n) =1$ for all $n \in {\Bbb N}$;
\item[(b)] $\Theta_{\ell_2}(n) = O(n^{-1/2}).$
\end{itemize}

In this paper we will show that $\Theta_{\ell_p}(n) \simeq n^{-1/p}$ and the same is true for an infinite $\ell_p$ sum of finite dimensional spaces. That result will be a consequence of the estimation of 
upper and lower bounds for covering index $\Theta_X$ in terms of the geometry of $X$, namely those bounds depend  
on the asymptotic moduli of uniform convexity and uniform smoothness for AUC or AUS renormings of $X$.
Moreover, we will show that for the Tsirelson space $T$, the index of its dual satisfies 
$\Theta_{T^*}(n) \geq 1/2$ for all $n \in {\Bbb N}$.\\

All the Banach spaces considered are real and separable. We will use the symbol 
$\lesssim$ instead of Landau's $O$, in order to use 
$\gtrsim$ for the reverse estimation, that is, we write $f(t) \gtrsim g(t)$ if there exists $c>0$ such that $f(t) \geq c \, g(t)$ for every $t$ in the domain. We will use $\simeq$ in case $\lesssim$ and $\gtrsim$ hold simultaneously.
With the exception of these details, our notation is totally standard and we address the reader to generic references such as \cite{LT, JL, banach} for any unexplained definition.

\section{Upper bounds for $ \Theta_X(n) $}

Consider  a Banach space $X$ with a Schauder finite-dimensional decomposition $(E_n)_{n=1}^\infty$, abbreviated FDD, see \cite{LT} for precise definitions. The following notions are relative to a fixed FDD of $X$ for $1 \leq p,q <\infty$:
\begin{itemize}
\item[(1)] we say that $X$ has a upper $p$ estimate if there is $C>0$ such that
$$ \Big\| \sum_{i=1}^k x_i \Big\|  \leq C \Big(  \sum_{i=1}^k \|x_i\|^p \Big)^{1/p} $$
for any $k \in {\Bbb N}$ and choice of $x_i \in E_i$ for $1 \leq i \leq k$;
\item[(2)] we say that $X$ has a lower $q$ estimate if there is $c>0$ such that
$$ \Big\| \sum_{i=1}^k x_i \Big\|  \geq c \Big(  \sum_{i=1}^k \|x_i\|^q \Big)^{1/q} $$
for any $k \in {\Bbb N}$ and choice of $x_i \in E_i$ for $1 \leq i \leq k$.
\end{itemize}

The idea behind the proof of the following result comes from \cite[Example 2.1]{raja2} were we used  ``parabolic cylinders'' to cut the unit Hilbert ball.

\begin{theo}\label{upper}
If $X$ has a shrinking FDD satisfying a lower $q$ estimate for some $q \geq 1$, then $ \Theta_X(n) \lesssim n^{-1/q} $.
\end{theo}

\noindent
\begin{proof}
Without loss of generality, we will suppose that the FDD indexed by ${\Bbb N} \cup \{0\}$, with $E_0$ of dimension $1$, so we can identify $E_0={\Bbb R}$. For every $x \in X$ we will write $(x_i)$ for its decomposition according to $(E_i)$.
Note that under our hypotheses, the FDD induces a bounded linear operator
$$ X \longrightarrow \Big(\bigoplus_{i=0}^\infty E_i \Big)_{\ell_q} . $$
We will estimate the covering index $\Theta_X$ using a particular decomposition into an even number of pieces $2n$.
Now, for $1 \leq j \leq 2n$, consider the set
$$ A_j = \Big\{ x \in B_X : 
(-1)^j x_0 \leq 1/2 - c^q n \sum_{m=0}^\infty \| x_{2nm+j} \|^q \Big\} ,$$
where $c>0$ is the constant associated to the lower $q$ estimate.
These sets are evidently closed and convex. Moreover, $\bigcup_{j=1}^{2n} A_j=B_X$. 
Indeed, assume that some $x =(x_i) \in B_X$  belongs to none of the sets. Then we have the reversed inequalities 
$$ (-1)^j x_0 > 1/2 - c^q n \sum_{n=0}^\infty \| x_{2kn+j} \|^q ,$$
whose sum over $1 \leq j \leq 2n$ gives
$$ 0 > n -c^q n \sum_{i=1}^\infty \| x_{i} \|^q .$$
We deduce
$$  1 <   \, c^q  \sum_{i=1}^\infty \| x_{i} \|^q \,  \leq \|x\|^q  $$
that is a contradiction with $x \in B_X$.
Now we will estimate the radius of a finite codimensional ball contained in $A_j$ for any $j$. 
Let $Y \subset X$ be a finite codimensional subspace. Since the FDD is shrinking, up to a small perturbation we may assume that  annihilator of $Y$ is contained in $\bigoplus_{i=0}^N E_i^*$. Let $x \in A_j$ be the center of the ball. We may take $m \in {\Bbb N}$ such that $2mn+j > N$ and the coordinate $x_{2mn+j }$ is negligible. Therefore, if $x + y \in A_j$, for the estimation of $\|y\|$ only the coordinates $0$ and $2mn+j $ are relevant. We have
$$ c^qn \|y_{2mn+j} \|^q \leq 1/2 - (-1)^{j}x_0 \leq 1/2 + |x_0| \leq M $$
where $M>0$ depends only on the norm of the projection onto $E_0$. We get that
$$ \| y_{2mn+j} \| \leq \Big( \frac{M}{c^q n} \Big)^{1/q} .$$
That estimation shows 
$ \varrho(A_j) \lesssim n^{-1/q} $
as wished.
\end{proof}

\begin{rema}\label{rem_sub}
We can extend Theorem \ref{upper} to Banach spaces having an overspace with a shrinking FDD satisfying a lower $q$ estimate, provided that the FDD also satisfy an upper $q$ estimate. Indeed, instead of reducing $y$ to just one coordinate we can work with the whole pack thanks to the upper estimation. 
\end{rema}

We will recall the definition of  the {\it modulus of asymptotic uniform convexity} of a Banach space $X$, notion coined in \cite{JLPS} but introduced by Milman  \cite{Milman} under a different name, see also \cite{LB}. The AUC modulus of $X$ is defined for 
$\varepsilon>0$ as
$$ \overline{\delta}_{X}(\varepsilon) =
\inf_{\|x\|=1} \sup_{{\mathcal D}(X/Y)<\infty} \inf_{y \in Y, \, \|y\| \geq \varepsilon} ( \|x+y\|-1 ). $$
The space is said to be {\it asymptotically uniformly convex} (AUC) if $\overline{\delta}_{X}(\varepsilon)>0$
for every $\varepsilon>0$. Finally, the space $X$ is said $q$-AUC for some $q \geq 1$ if 
$\overline{\delta}_{X}(\varepsilon)  \gtrsim \varepsilon^q$.
According to \cite{KOS} for every separable reflexive AUC space there exists a renorming making it $q$-AUC.

\begin{coro}
Let $X$ be a reflexive Banach space with a FDD. If $X$ admits an equivalent $q$-AUC norm, then 
$ \Theta_X(n) \lesssim n^{-1/q} $.
\end{coro}

\noindent
\begin{proof}
According to a result of Prus \cite{Prus} (see also \cite{JLPS, LB}), there exists a blocking of the FDD, that is a further FDD, satisfying a lower $q$ estimate.
\end{proof}

\section{Lower bounds for $ \Theta_X(n) $}

The main tool in this section will be the {\it goal derivation} studied in \cite{raja1}.
For a set $A \in {\mathcal C}(X)$ and $\varepsilon>0$ we define 
$$ [ A ]'_\varepsilon=\{x \in A: \forall \, U \, w \mbox{-neighbourhood of~} x,
\varrho(A \cap U) > \varepsilon \}.$$
It is not evident but true that $ [ A ]'_\varepsilon \in  {\mathcal C}(X)$ too.
We can define derived sets of superior order by taking $ [ A ]^n_{\varepsilon} = [ [A]^{n-1}_{\varepsilon} ]'_{\varepsilon} $, 
for any $n \in {\Bbb N}$ and eventually for ordinal numbers.
The {\it goal Szlenk index}, denoted $Gz(A,\varepsilon)$ is defined as
$$ Gz(A,\varepsilon) = \inf\{n \in {\Bbb N}: [ A ]^n_{\varepsilon} =\emptyset \} $$
when it is finite. Otherwise, we will simply say that $Gz(A,\varepsilon) \geq \omega$ (here we are not interested in dealing with ordinals). Actually, the original definition of the goal derivation uses ``$\geq$'' instead of ``$>$'' as above, but that does not change the behaviour of $Gz(A,\varepsilon)$ as a function of $\varepsilon$.

\begin{prop}
$Gz(A,\varepsilon) \leq n$ if and only if there is a covering 
$(A_k)_{k=1}^n$ of $B_X$ by weakly closed sets such that for every $x \in B_X$ there is $1 \leq k \leq n$ such that 
$$ \sup \{r>0: \, \exists Y \subset X, \, {\mathcal D}(X/Y) < \infty, ~x+rB_Y \subset A_k \} \leq \varepsilon. $$
\end{prop}

\noindent
\begin{proof} 
If $Gz(A,\varepsilon) \leq n$ just take $A_k=  [ A ]^{k-1}_\varepsilon$ for $1 \leq k \leq n$.
For the other implication, define the sets
$$ B_k = \big\{ x \in B_X: \big| \{ j: x \in A_j \} \big| > k \big\} $$
for $0 \leq k \leq n$ (obviously, $B_n =\emptyset$).
We claim that $[ A ]^k_\varepsilon \subset B_k$, that clearly implies the statement. 
The proof will be by induction. Evidently $B_0=B_X$. Assume now $k < n$ and
$[ A ]^k_\varepsilon \subset B_k$. Take any $x \in [ A ]^k_\varepsilon \setminus B_{k+1}$. Then there is a subset $I \subset \{1,\dots,n\}$ with $|I|=k+1$ such that $x \in A_i $ if and only if $i \in I$. That implies
$$ V = [ A ]^k_\varepsilon \setminus \bigcup_{j \not \in I} A_j \subset \bigcap_{i \in I} A_i. $$
By construction $V$ is a relatively weakly open neighbourhood of $x$ such that 
$$ \sup \{r>0: \, \exists Y \subset X, \, {\mathcal D}(X/Y) < \infty, ~x+rB_Y \subset V \} < \varepsilon $$
since the hypothesis is fulfilled for some $i \in I$. But the same is true for any other point in $V$, meaning
$\varrho(V) \leq \varepsilon$. That implies $ [ A ]^{k+1}_\varepsilon \subset B_{k+1}$ as wished.
\end{proof}

\begin{coro}\label{coro_gz}
If $\varepsilon > \Theta_X(n)$, then $Gz(B_X, \varepsilon) \leq n$.
\end{coro}

Since $Gz(B_X,\varepsilon) \gtrsim n^{-1}$, see \cite{raja1}, that gives a first lower bound for $\Theta_X$.

\begin{coro}
For any Banach space $X$, we have $\Theta_X(n) \gtrsim n^{-1}$.
\end{coro}

In order to give more precise lower bound we need the notion of {\it modulus of asymptotic uniform  smoothness} of a Banach space. As the AUC modulus, the AUS modulus was coined in \cite{JLPS} but introduced by Milman  \cite{Milman} with a different name, see also \cite{LB}. The AUC modulus of $X$ is defined for 
$\varepsilon>0$ as
$$ \overline{\rho}_{X}(\varepsilon) =
\sup_{\|x\|=1} \inf_{{\mathcal D}(X/Y)<\infty} \sup_{y \in Y, \, \|y\| \leq \varepsilon} ( \|x+y\|-1 ). $$
The space $X$ is said to be
{\it asymptotically uniformly smooth} if $\lim_{\varepsilon \rightarrow 0}
\varepsilon^{-1}\overline{\rho}_{X}(\varepsilon) =0$. Finally, the space $X$ is said $p$-AUS for some $p > 1$ 
if $\overline{\rho}_{X}(\varepsilon) \lesssim \varepsilon^q$.
According to \cite{KOS} for every separable AUS space there exists a renorming making it $p$-AUS. Although we do not have AUS property for $p=1$ we will consider the notion of ``$1$-AUS'' for the next result from \cite{raja1}.

\begin{prop} 
Suppose that $X$ is $p$-AUS for $p \geq 1$, then $Gz(B_X,\varepsilon) \gtrsim \varepsilon^{-p}$.
\end{prop}

The straightforward combination of Corollary \ref{coro_gz} and the previous proposition gives the following lower bound for 
$\Theta_X$.

\begin{theo}\label{lower}
If $X$ is $p$-AUS renormable for $p \geq 1$, then $ \Theta_X(n) \gtrsim n^{-1/p} $.
\end{theo}

\section{Asymptotic estimations and problems}

The combination of upper and lower bounds will allow us to give precise asymptotic estimations of $\Theta_X$. 
Since we are interested in isomorphic results, the notion of equivalence for sequences fits to our purpose.
We say that the sequences 
$(a_n), (b_n) \subset {\Bbb R}^+$ are equivalent, denoted as $a_n \simeq b_n$, if there exist constants $c,d>0$ such that
$c \, a_n \leq b_n \leq d\, a_n$ for all $n \in {\Bbb N}$. 

\begin{prop}
Let $X$ be isomorphic to a subspace of an $\ell_p$ sum of finite-dimensional normed spaces. 
Then $\Theta_X(n) \simeq n^{-1/p}$.
\end{prop}

\noindent
\begin{proof} 
This is just  a combination of Theorem \ref{upper}, Remark \ref{rem_sub} and Theorem \ref{lower}.
\end{proof}

In particular, the $\Theta_X$ index distinguishes between $\ell_p$ spaces (in particular, better than the type or the cotype).

\begin{coro}
$\Theta_{\ell_p} (n) \simeq n^{-1/p}$ for $1 \leq p <+\infty$.
\end{coro}

A deep result from \cite{JLPS} says that a reflexive Banach space that admits both $p$-AUC and $p$-AUS renormings can be embedded as a subspace of an $\ell_p$ sum of finite-dimensional normed spaces. 

\begin{coro}
Let $X$ be a reflexive Banach space and $p>1$ such that $X$ admits  $p$-AUC and $p$-AUS renormings, then 
$\Theta_{X} (n) \simeq n^{-1/p}$.
\end{coro}

Recall that for a general Banach space $(\Theta_X(n))_{n=1}^\infty$ may not converge to $0$.

\begin{prop}
Let $\varepsilon \in (0,1]$ be such that $Gz(B_{X},\eta) \geq \omega$ for every $\eta \in (0,\varepsilon)$.
Then $\lim_n \Theta_{X}(n) \geq \varepsilon$.
\end{prop}

\noindent
\begin{proof}
Indeed, if $\eta <\varepsilon$, then $\Theta_{X}(n) \geq \eta$ for all $n \in {\Bbb N}$ by Corollary \ref{coro_gz}.
\end{proof}

We retrieve the result mentioned in the introduction that $\Theta_{c_0}(n)=1$ for all $n \in {\Bbb N}$, since the goal set derivation satisfies $[B_{c_0}]'_\eta = B_{c_0}$ for every $\eta \in (0,1)$ so $Gz(B_{c_0},\eta)$ is not even defined (we put 
``$\infty$'' as something above all the ordinals).\\

The following example shows that the behaviour of $\Theta_X$ can be very diverse among reflexive Banach spaces. Consider the dual $T^*$ of the Tsirelson space $T$, some times also called Tsirelson space, see [8] for a construction. In order to avoid confusion, let us say that $T$ looks like $\ell_1$ and $T^*$ takes after $c_0$.

\begin{coro}
Let $T$ be the Tsirelson space, then $\lim_n \Theta_{T^*}(n) \geq 1/2$.
\end{coro}

\noindent
\begin{proof}
Indeed, if $\varepsilon <1/2$, then $Gz(B_{T^*},\varepsilon) \geq \omega$ by \cite[Theorem 4.4]{raja1}.
\end{proof}

I would like to finish with some general comments and open problems. 
The index $\Theta_X$ is defined in a very simple way, but it can recognize structural features of Banach spaces. 
Our initial reason for estimating the radius of finite codimensional balls is the fact that it essentially ``survive'' through uniform homeomorphisms. That is the so called the {\it Gorelik principle}, see \cite{LB, GKL}, that was used by Godefroy, Kalton and Lancien \cite{GKL} to prove that AUS renormabilty  of a Banach space is a uniform invariant. 
However, their proof makes a heavy use of duality so it would be desirable a better understanding in terms of the space itself.
Despite this motivation, we feel that our results are ``too much linear'' to combine well with Gorelik's.\\ 

There is an evident parallelism of $\Theta_X$ with the {\it entropy numbers}. Given a relatively compact subset $K \subset X$, the entropy number $ {\mathcal E}_n(K)$ is defined as the infimum of the $\varepsilon>0$ such that $K$ can be covered with $n$ balls of radius $\varepsilon$. Not just parallelism but a bit of interaction: Let $T: X \rightarrow Z$ be a compact linear operator such that $\lim_n {\mathcal E}_n(T(B_X))/\Theta_X(n)=0$. It is not difficult to prove that for every $\varepsilon>0$ there is finite-codimensional subspace $Y \subset X$ such that $\|T|_{Y}\| < \varepsilon$, meaning that such $T$ can be approximated by finite rank operators.\\

Now, we will explicitly formulate three problems that we could not solve.

\begin{open}
Find the exact value of $\Theta_{\ell_2}(2)$.
\end{open}

We know that $0.707 \leq \Theta_{\ell_2}(2) \leq 0.931$, see \cite{raja2} for the upper bound. The lower bound can be obtained with a suitable choice of a $1$-codimensional ball (just make a $2$-dimensional picture).

\begin{open}
Does  $\Theta_X$ really depends on both the moduli of AUC and AUS?
\end{open}

Maybe there is a connection to the $\overline{\beta}_X$ modulus of Rolewicz, that takes an intermediate place between the AUC and AUS moduli, see \cite{DKLR, LB}. 
\begin{open}
If $\lim_n \Theta_X(n) >0$, is $c_0$ crudely finitely representable in $X$?
\end{open}

The answer is affirmative if $Gz(B_X,\varepsilon) \geq \omega$ for some $\varepsilon>0$, see \cite{raja1}.

\vspace{1cm}

\begin{flushright}
Departamento de Matem\'aticas\\ Universidad de Murcia\\
Campus de Espinardo\\ 30100 Espinardo, Murcia, SPAIN\\
E-mail: matias@um.es
\end{flushright}

\end{document}